%% file: Note-C-surgery.tex
\documentclass[10pt]{amsart}
\usepackage{hyperref,amsmath,amsthm,amsfonts,amscd,flafter,epsf,amssymb,diagrams}
\usepackage{epsfig}
\usepackage{amsfonts}
\usepackage{amssymb}
\usepackage{amsmath}
\usepackage{amsthm}
\usepackage{latexsym}
\usepackage{amscd}
\usepackage{flafter}
\usepackage{epsf}
\input{diagrams.sty}
 \input{command1}
\begin{document}

\title[A note on splicing formulas]{ A note on combinatorial splicing formulas for Heegaard Floer homology}%
\author{Eaman Eftekhary}%
\address{School of Mathematics, Institute for Studies in Theoretical Physics and Mathematics (IPM),
P. O. Box 19395-5746, Tehran, Iran}%
\email{eaman@ipm.ir}
%\thanks{The authors is partially supported by a NSF grant}%
%\subjclass{}%

%\date{December 2005}%
%\dedicatory{}%
%\commby{}%
\maketitle
% ----------------------------------------------------------------
\begin{abstract}
We give a precise description of splicing formulas from a previous paper in terms
of knot Flor complex associated with a knot in a homology sphere.
\end{abstract}
\section{Introduction}
Heegaard Floer homology has provided topologists with a quite powerful technique in the study of
surgery along null-homologous knots inside three-manifolds. After the introduction of surgery formulas
for Heegaard Floer homology by Ozsv\'ath and Szab\'o (\cite{OS-surgery,OS-Qsurgery}),
several interesting results have been proved (e.g. \cite{OS-Qsurgery,Hed}).
The author has considered another type of surgery formulas which are suggested from
the combinatorial approach of Sarkar and Wang to Heegaard Floer homology (see \cite{Ef-C-surgery}, also \cite{SW}, also
\cite{MOS,MOST}). 
More generally, there is a formula for Heegaard Floer homology of the three-manifold 
obtained by splicing two knot complements, provided in \cite{Ef-splicing}. 
In this paper, we will compare these formulas, and will bring the
combinatorial surgery formulas into a form similar to the formulas of Ozsv\'ath and Szab\'o.
More interestingly, this will give a more explicit formula for Heegaard Floer homology of the
three-manifold obtained by splicing two knot complements, which is described in terms of
knot Floer complexes of the two knots.\\

\begin{figure}
\mbox{\vbox{\epsfbox{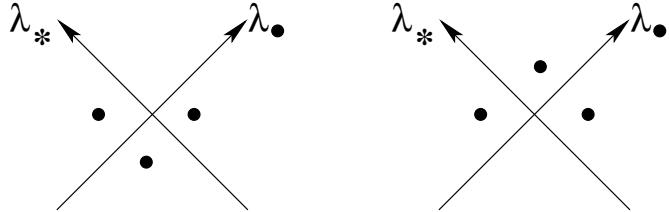}}} \caption{\label{fig:punctures}
{For defining chain maps between $C_\bullet(K)$ and $C_\star(K)$, the
punctures around the intersection point of $\lambda_\bullet$ and $\lambda_\star$
should be chosen as illustrated in the above diagrams.
}}
\end{figure}

Let $K$ be a knot inside
the homology sphere $Y$. We may remove a tubular neighborhood of $K$ and
glue it back to obtain the three-manifold $Y_{p/q}=Y_{p/q}(K)$, which is the
result of $p/q$-surgery on $K$. The core of the solid torus, which is the
tubular neighborhood of $K$, will represent a knot in $Y_{p/q}$ which will be
denoted by $K_{p/q}$. We may denote $(Y,K)$ by $(Y_\infty,K_\infty)$, as
an extension of the above notation. Let $\Hbb_\bullet(K)$ be the Heegaard
Floer homology group $\widehat{\text{HFK}}(Y_\bullet,K_\bullet)$ for
$\bullet\in \mathbb{Q}\cup\{\infty\}$. Note that $\widehat{\text{HFK}}$ is
defined for knots inside rational homology spheres (see \cite{OS-Qsurgery}), and that
$\Hbb_0(K)=\widehat{\text{HFL}}(Y,K)$ is the longitude Floer homology
of $K$ from \cite{Ef-Whitehead}. In all these cases, we choose the
coefficient ring to be $\Z/2\Z$ (for simplicity).
If we choose a Heegaard diagram for $Y-K$ and let $\lambda_\bullet$ denote
a longitude which has framing $\bullet\in \Z\cup\{\infty\}$ (with $\lambda_\infty=\mu$
the meridian for $K$), one may observe that the pairs $(\lambda_\infty,\lambda_1)$ and
$(\lambda_1,\lambda_0)$ have a single intersection point in the Heegaard diagram.
Let $(\bullet,\star)\in\{(\infty,1),(1,0)\}$ correspond to either of these pairs.
There are four quadrants around the intersection point of $\lambda_\star$ and
$\lambda_\bullet$. If we puncture three of these quadrants and consider the
corresponding holomorphic triangle map, we obtain an induced map $\Hbb_\bullet\ra
\Hbb_\star$. If the punctures are chosen as in figure~\ref{fig:punctures}, the
result would be  two maps
$\phi,\phibar:\Hbb_\infty(K)\ra \Hbb_1(K)$ and two other maps
$\psi,\psibar:\Hbb_1(K)\ra \Hbb_0(K)$ so that the following
two sequences are exact:
\begin{displaymath}
\begin{split}
&\Hbb_0(K) \xra{\psi}\Hbb_1(K)\xra{\phibar}\Hbb_\infty(K),\ \ \ \text{and}\\
&\Hbb_0(K) \xra{\psibar}\Hbb_1(K)\xra{\phi}\Hbb_\infty(K).\\
\end{split}
\end{displaymath}
The homology of the mapping cones of $\phi$ (or $\phibar$) and $\psi$ (or $\psibar$)
are $\Hbb_0(K)$ and $\Hbb_\infty(K)$ respectively (see \cite{Ef-splicing}).
With the above notation fixed, we have proved the following surgery formula in \cite{Ef-splicing} (The
orientation convention of \cite{Ef-splicing} is different from that of \cite{Ef-surgery} and \cite{OS-surgery}, so we have
changed the direction of all maps to compensate this difference):
\begin{thm}
 Suppose that $K_1$ and $K_2$ are two knots inside homology spheres $Y_1$ and $Y_2$ respectively.
Define the maps $\phi^i,\phibar^i,\psi^i$ and $\psibar^i$ for $K_i$ as above.
Let $\M=\M(K_1,K_2)$ denote the following cube of maps
\begin{displaymath}
\begin{diagram}
\Hinfinf   &              & \lTo^{I\otimes\phi^2}&&             & \Hinfone       &               &          \\
	   &\luTo^{\phi^1\otimes I}&              &&             & \uTo         &\luTo^{\phi^1\otimes I} &          \\
\uTo^{\Gtt}&              & \Honeinf     &&\lTo^{I\otimes\phi^2}& \HonV          &               &\Honeone  \\
           &              &              &&             & \vLine_{\Got}    &               &          \\
	   &              & \uTo^{\Gto}  &&             &                &               &          \\
\Hzerozero & \hLine         & \VonH       &\rTo{\psi^1\otimes I}&& \Honezero      &               &\dTo^{Id} \\
	  & \rdTo_{I\otimes \psi^2}&              &&             &                &\rdTo_{I\otimes \psi^2} &          \\
	   &              &\Hzeroone     &&             &\rTo^{\psi^1\otimes I}   &               &\Honeone , \\
\end{diagram}
\end{displaymath}
where $\Hbb_{\star}^1\otimes\Hbb_{\bullet}^2=\Hbb_\star(K_1)\otimes_{\Z/2\Z}\Hbb_\bullet(K_2)$ for $\star,\bullet\in\{\infty,1,0\}$.
The differential $d_{\M}$ of the complex $\M$ is defined to be the sum of all the maps that appear in
this cube.
Then the Heegaard Floer homology of the three-manifold $Y$, obtained by splicing knot complements
$Y_1-K_1$ and $Y_2-K_2$, is given by
$$\widehat{\text{\emph{HF}}}(Y;\Z/2\Z)=H_*(\M,d_{\M}),$$
where $H_*(\M,d_{\M})$ denotes the homology of the cube $\M$.
\end{thm}

On the other hand, the groups $\Hbb_\bullet (K_i)$ for $\bullet\in \Z$ are described in terms of knot
Floer complex for $(Y_i,K_i)$ in \cite{Ef-surgery}. We have shown that the following is true:
\begin{thm}
Suppose that $(Y,K)$ is a homology sphere and ${\B}$ is the complex
$\widehat{\text{\emph{CF}}}(Y,\Z/2\Z)$, equipped by the filtration induced by $K$.
 Let ${\B}(\geq s)$ denote the
quotient-complex generated by those generators whose relative $\SpinC$
class in $\Z=\RelSpinC(Y,K)$ is greater than or equal to $s\in\Z$. The homology
group $\widehat{\text{\emph{HFK}}}(Y_n(K),K_n;s;\Z/2\Z)$ are isomorphic to
the homology of the following chain complex
$${\B}(>s-n)\xra{\imath}\B\xleftarrow{\imath}\B(\geq -s),$$
where $\imath$ is the inclusion map.
\end{thm}

The question we would like to address in this paper is the description of the maps
$\phi,\phibar,\psi$ and $\psibar$ in terms of a presentation of $\Hbb_\bullet(K)$
given by the above theorem. This would give surgery and splicing formulas
which will look similar to the surgery formulas of Ozsv\'ath and Szab\'o in \cite{OS-surgery,OS-Qsurgery}.\\

Namely, if we denote the part of $\Hbb_\bullet(K)$ in relative $\SpinC$ class
$s\in \Z=\RelSpinC(Y,K)$ by $\Hbb_\bullet(K,s)$ we will prove
\begin{thm}\label{thm:main}
Under the identification of $\Hbb_\infty(K,s)$ with $\B\{s\}$, $\Hbb_1(K,s)$ with the
homology of the complex
$$C_1(s)=\Big({\B}\{\geq s\}\xra{\imath}\B\xleftarrow{\imath}\B\{\geq -s\}\Big),$$
and $\Hbb_0(K,s)$ with the homology of the complex
$$C_0(s)=\Big({\B}\{\geq s+1\}\xra{\imath}\B\xleftarrow{\imath}\B\{\geq -s\}\Big),$$
the map $\phi:\Hbb_1(K)\ra \Hbb_\infty(K)$ is induced by the sum of
maps $\phi_s$ which take the quotient complex
$\B\{\geq s\}$ of the complex $C_1(s)$ to it quotient $\B\{s\}=\Hbb_\infty(K,s)$. The map
$\phibar$ is the sum of maps $\phibar_s$ which are induced by using
the map $\phi_{-s}$ followed by the isomorphism $\B\{s\}\simeq \B\{-s\}$.
Moreover, the map $\psi$ is a sum of inclusion maps
$\psi_s$ from $C_0(s-1)$ to $C_1(s)$ and the
 map $\psibar$ is induced by the sum of inclusion maps $\psibar_s$ from $C_0(s)$ to $C_1(s)$.
\end{thm}

\section{Surgery on null-homologous knots; A short review}
In this section we will review the construction of \cite{Ef-surgery}, at least to the extent relevant
for the purposes of this paper. The reader is referred to \cite{Ef-surgery} for the proofs.\\

Consider a Heegaard diagram for the pair $(Y,K)$. Suppose that
the curve $\mu=\beta_g$ in the Heegaard diagram
$$H=(\Sig,\alphas=\{\alpha_1,...,\alpha_g\},\betas=\{\beta_1,...,
\beta_g\},p)$$
corresponds to the meridian of $K$ and that the marked point $p$ is
placed on $\beta_g$. One may assume that the curve $\beta_g$ cuts
$\alpha_g$ once and that this is the only element of $\alphas$ that has
an intersection point with $\beta_g$. Suppose that $\lambda$
represents a longitude for the knot $K$ (i.e. it cuts $\beta_g$ once
and stays disjoint from other elements of $\betas$) such that the
Heegaard diagram
$$(\Sig,\alphas,\{\beta_1,...,\beta_{g-1},\lambda\})$$
represents the three-manifold $Y_0(K)$ obtained by zero surgery on $K$. Winding $\lambda$ around
$\beta_g$- if it is done $n$ times- would produce a Heegaard diagram
for the three-manifold $Y_n(K)$. More precisely, if the resulting
curve is denoted by $\lambda_n$, the Heegaard diagram
$$H_n=(\Sig,\alphas,\betas_n=\{\beta_1,...,\beta_{g-1},\lambda_n\},p_n)$$
would give a diagram associated with the knot $(Y_n(K),K_n)$, where
$y=p_n$ is a marked point at the intersection of $\lambda_n$ and $\beta_g$.\\

The curve $\lambda_n$ intersects the $\alpha$-curve $\alpha_g$ in
$n$-points which appear in the winding region (there may be other
intersections outside the winding region). Denote these points
of intersection by
$$...,x_{-2},x_{-1},x_{0},x_{1},x_{2},...,$$
where $x_1$ is the intersection point with the property that three of
its four neighboring quadrants belong to the regions that
 contain either $p_n$ as a corner . Any
generator which is supported in the winding region is of the form
$$\{x_i\}\cup \y_0=\{y_1,...,y_{g-1},x_i\},$$
and it is in correspondence with the generator
$$\y=\{x\}\cup\y_0=\{y_1,...,y_{g-1},x\}$$
for the complex associated with the knot $(Y,K)$, where $x$ denotes the
unique intersection point of $\alpha_g$ and $\beta_g$. Denote the former
generator by $(\y)_i$, keeping track of the intersection point $x_i$
among those in the winding region.\\

By taking an intersection point $\x$ in relative $\SpinC$ class $\relspinc-k\in \RelSpinC(Y,K)=\Z$ to $(\x)_k$  when $k\leq 0$ and $(\x)_{k+n}$ if $k>0$ we obtain 
an identification of the complex $\widehat{\CFKT}(Y_n(K),K_n;\relspinc)$ and the complex
$$\ov{\U}(\relspinc)=\B\{\geq \relspinc\}\oplus \B\{>n-\relspinc\}=\B\{\geq \relspinc\}$$
for large values of $n$, having fixed $\relspinc\in \Z$.

%%%%%%%%%%%%%%%%%%%%%%%%%%
%%%%%%%%%%%%%%%%%%%%%%%%%%
%%%%%%%%%%%%%%%%%%%%%%%%%%
%%%%%%%%%%%%%%%%%%%%%%%%%%
%%%%%%%%%%%%%%%%%%%%%%%%%%
%%%%%%%%%%%%%%%%%%%%%%%%%%
%%%%%%%%%%%%%%%%%%%%%%%%%%
%%%%%%%%%%%%%%%%%%%%%%%%%%
%%%%%%%%%%%%%%%%%%%%%%%%%%
%%%%%%%%%%%%%%%%%%%%%%%%%%
%%%%%%%%%%%%%%%%%%%%%%%%%%
%%%%%%%%%%%%%%%%%%%%%%%%%%
%%%%%%%%%%%%%%%%%%%%%%%%%%
%%%%%%%%%%%%%%%%%%%%%%%%%%
%%%%%%%%%%%%%%%%%%%%%%%%%%
%%%%%%%%%%%%%%%%%%%%%%%%%%

\begin{figure}
\mbox{\vbox{\epsfbox{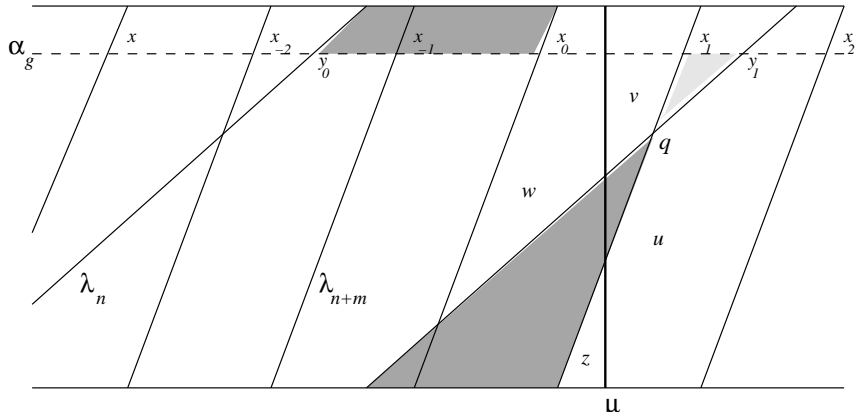}}} \caption{\label{fig:HD-Rmn} {
The Heegaard diagram $R_{m,n}$. The shaded triangles are $\Delta_0$
and $\Delta_1$. The marked points $u,v,w$ and $z$ are placed in
$D_1^R,D_2^R,D_2^L$ and $D_1^L$ respectively.
}}
\end{figure}

 Choose an intersection point between the curves $\lambda_n$
and $\lambda_{m+n}$ in the middle of the winding region, denoted by
$q$. We will assume that $m=\ell n$ for an integer $\ell$ which is chosen to
be appropriately large. We continue to assume that $\alpha_g$ is the
unique $\alpha$-curve in the winding region.
From the $4$ quadrants around the intersection
point $q$, two of them are parts of  small triangles $\Delta_0$ and
$\Delta_1$ between
$\alphas,\betas_n$ and $\betas_{m+n}$.
We may assume that the intersection points between $\alpha_g$ and
$\lambda_{m+n}$ in the winding region are
$$...,y_{-2},y_{-1},y_0,y_1,y_2,...,$$
and that the intersection points between $\lambda_n$ and $\alpha_g$
are $$...,x_{-2},x_{-1},x_0,x_1,x_2,...$$
as before. We may also assume that the domain
$\Delta_i$ for $i=0,1$ is the triangle with vertices
$q,x_i$ and $y_i$, and $\Delta_1$ is one of the connected
domains in the complement of curves $\Sig\setminus C$ where
 $$C=\alphas \cup\betas_n\cup \betas_{m+n}.$$
Other that $\Delta_0$ and $\Delta_1$ there are two other
domains which have $q$ as a corner. One of them is on the
right-hand-side of both $\lambda_n$ and $\lambda_{m+n}$, denoted by
$D_1$, and the other one is on the left-hand-side of both of them,
denoted by $D_2$. The domains
$D_1$ and $D_2$ are assumed to be connected regions
in the complement of the
curves $\Sig\setminus C$.
We may assume that the meridian $\mu$ passes through the regions
$D_1,D_2$ and $\Delta_0$, cutting each of them into two parts:
$\Delta_0=\Delta_0^R\cup\Delta_0^L$, $D_1=D_1^R\cup D_1^L$ and
$D_2=D_2^R\cup D_2^L$. Here $\Delta_0^R\subset \Delta_0$ is the
part on the right-hand-side of $\mu$ and $\Delta_0^L$ is the part
on the left-hand-side. Similarly for the other partitions. Choose
the marked points so that $u$ is in $D_1^R$, $v$ is in $D_2^R$, $w$
is in $D_2^L$ and $z$ is in $D_1^L$ (see figure~\ref{fig:HD-Rmn}).
We obtain the Heegaard diagram
$$R_{m,n}=(\Sig,\alphas,\betas,\betas_n,\betas_{m+n};u,v,w,z),$$
which determines three chain complexes associated with the pairs
$(\alphas,\betas),(\alphas,\betas_n),$ and $(\alphas,\betas_{m+n})$
once we use the marked points $u,v,w,z$. The complexes associated
with the last two pairs are simply $\widehat{\CFKT}(Y_n(K),K_n)$ and
$\widehat{\CFKT}(Y_{m+n}(K),K_{m+n})$, obtained by puncturing the
surface at these marked points. The complex associated with the first pair
is constructed by making punctures at $u,v$, counting holomorphic disks,
and twisting according to the intersection number of the disks with either
of the marked points $w,z$. 
 As we have
already seen, $\ov{\CFKT}(K_{m+n};\relspinc)$ may be identified with the complex
 $\ov{\U}(\relspinc)$ if $m$ is large enough. In fact,  the holomorphic triangle map gives a 
 long sequence for each relative $\SpinC$ class $\relspinc\in \Z$ as 
 \begin{displaymath}
 \begin{diagram}
  \hdots&\rTo{f^{\relspinc}}&\widehat{\text{CF}}(Y)
&\rTo{h^\relspinc}&\widehat{\CFKT}(K_n,\relspinc)&\rTo{g^\relspinc}
&\ov{\E}(\relspinc)&\rTo{f^{\relspinc}}&
\widehat{\text{CF}}(Y)&\rTo{h}&\hdots,
 \end{diagram}
\end{displaymath}
such that the induced maps in homology form a long exact sequence, if $m=\ell.n$ and $\ell$ is large enough.
Here $\ov{\E}(\relspinc)$ is the complex $\B\{\geq \relspinc\}\oplus \B\{>n-\relspinc\}$. This
implies that $\ov{\text{HFK}}(Y_n(K),K_n;\relspinc)$ may be computed as the
homology of the mapping cone of $f^{\relspinc}:\ov{\E}(\relspinc)\ra \B$,
which is the sum of two inclusion maps.
One particular case, is the case where $n=1$ and this later mapping cone is the same as the complex
$C_1(\relspinc)$ mentioned before.\\

%%%%%%%%%%%%%%%%%%%%%%%%%%
%%%%%%%%%%%%%%%%%%%%%%%%%%
%%%%%%%%%%%%%%%%%%%%%%%%%%
%%%%%%%%%%%%%%%%%%%%%%%%%%
%%%%%%%%%%%%%%%%%%%%%%%%%%
%%%%%%%%%%%%%%%%%%%%%%%%%%
%%%%%%%%%%%%%%%%%%%%%%%%%%
%%%%%%%%%%%%%%%%%%%%%%%%%%
%%%%%%%%%%%%%%%%%%%%%%%%%%
%%%%%%%%%%%%%%%%%%%%%%%%%%
%%%%%%%%%%%%%%%%%%%%%%%%%%
%%%%%%%%%%%%%%%%%%%%%%%%%%
%%%%%%%%%%%%%%%%%%%%%%%%%%
%%%%%%%%%%%%%%%%%%%%%%%%%%
%%%%%%%%%%%%%%%%%%%%%%%%%%
%%%%%%%%%%%%%%%%%%%%%%%%%%

\section{Understanding the maps}
In this section we prove theorem~\ref{thm:main}.\\
\begin{proof}
In the construction of the previous section,
note that $A(\relspinc)=\ov{\CFKT}(Y_n(K),K_n;\relspinc)$ may be written as the mapping cone of
some chain map $\rho:A_1(\relspinc)\ra A_2(\relspinc)$, where $A_1(\relspinc)$ is the part of $A(\relspinc)$ generated
by the generators of the form $(\x)_1$. The differential of $A(\relspinc)$ induces a differential on the quotient
complex $A_1(\relspinc)$ and the sub-complex $A_2(\relspinc)$ of it. There is a similar decomposition of
$B=\ov{\CFKT}(Y_{m+n}(K),K_{m+n})$ as the mapping cone of $\sigma:B_1\ra B_2$, and the image of $A_1$ under $g$
is in $B_1$. Under the identification of $B(\relspinc)$ with $\ov{\U}(\relspinc)$, the quotient complex
$B_1(\relspinc)$ corresponds to the quotient complex $\B\{\relspinc\}$ of
$\B\{\geq \relspinc\}= \ov{\U}(\relspinc)$.\\

We will use the diagram $S_1$, together with the marked points $u,w$ and $z$
to define the map $\phi:C_1(s)=\widehat{\text{CFK}}(Y_1(K),K_1,s)\ra \widehat{\text{CFK}}(Y,K,s)=C_\infty(s)$.
Note that puncturing the surface $\Sig$ at these three marked points implies that this
chain map will respect the relative $\SpinC$ classes on the two sides. One can easily see
from the Heegaard diagram that the map $\phi$ is defined on the generators by
$\phi(\y)=
\x$, if $\y=(\x)_1+b_2$ for some generator $\x$ of $C_\infty$ and some $b_2\in B_2$.
This means that $\phi$ is trivial on the sub-complex $A_2$ of $A=C_1$.
The image of the generators  $(\x)_1$ of $A_1$ under the map $g$ is the generator $(\x)_1$, as a generator for
$\widehat{\CFKT}(Y_{m+1}(K),K_{m+1},s)$, i.e. $g$ identifies $A_1(\relspinc)$ with $B_1(\relspinc)$, via
the natural isomorphism of both of them with $\B\{\relspinc\}$.\\

On the other hand, since the map $h$ uses punctures at $u$ and $v$, its image does not contain
any generator of the form $(\x)_1$ in $\widehat{\CFKT}(Y_1(K),K_1,s)$, i.e. the image
of $h$ is in the sub-complex $A_2$ of $A$. This implies that if
an element $a$ in Ker$(g_*)$=Im$(h_*)$ is of the form $a_1\oplus a_2$ in $A=A_1\oplus A_2$, and 
$a_1=(\x)_1$, we should have
$(\x)_1=d(\y)_1=(dy)_1$ in $A_2$ for some generator $\y$ of $C_\infty$, i.e. $\x$ is trivial in 
$\Hbb_\infty(K)$.\\

Finally, note that by the exactness of the long sequence in homology we have
\begin{displaymath}
\begin{split}
\ov{\text{HFK}}(K_1,\relspinc)&=\text{Ker}(g^s_*)\oplus \text{Im}(g^s_*)\\
&=\text{Im}(h^s_*)\oplus {\text{Ker}}(f^{\relspinc}_*)\\
&=\frac{H_*(\B)}{\text{Im}(f^{\relspinc}_*)}\oplus  {\text{Ker}}(f^{\relspinc}_*)\oplus H_*(\frac{\ov{\E}(\relspinc)}{\text{Ker}(f^{\relspinc})}\xra{\simeq}\text{Im}(f^{\relspinc}))\\
&=H_*(\frac{\B}{\text{Im}(f^{\relspinc})}\xra{d_*}\text{Im}(f^{\relspinc})\xleftarrow{\simeq}
\frac{\ov{\E}(\relspinc)}{\text{Ker}(f^{\relspinc})}\xra{d_*} {\text{Ker}}(f^{\relspinc}))\\
&=H_*(\ov{\E}(\relspinc)\xra{f^{\relspinc}}\B).\\
\end{split}
\end{displaymath}
In these equalities, $H_*(\bullet)$ denotes the homology of the complex $\bullet$.
The above discussion shows that Ker$(g_*)$ is in the kernel of the induced map $\phi$. We thus need
to understand the image of Ker$(f^{\relspinc}_*)$ under the map $\phi$ in the
above description, which are represented in the last line of the above equality by the closed elements $e$
in $\ov{\E}(\relspinc)$ such that $f^{\relspinc}(e)$ is exact in $\B$. Any such element is the sum of
an exact element in $\ov{\E}(\relspinc)$ and an element of the form $g(a_1\oplus a_2)$ in $B_1$,
where $a_1\in A_1$ is closed
and $\rho(a_1)=da_2$ for some $a_2\in A_2$. Under the identification of $B_1$ with $A_1$, this image is of the
form $a_1\oplus e_2\in B_1\oplus B_2$. The discussion implies that for a closed element $e\oplus b\in \ov{\E}\oplus \B$
(i.e. such that $e$ is closed and $f(e)=db$) such that $e=a_1\oplus e_2\in B_1\oplus B_2$ we have
$\phi_*(e\oplus b)=a_1\in \widehat{\text{HFK}}(Y,K)=\Hbb_\infty(K)$. This proves the first assertion in theorem~\ref{thm:main}. The other claims follow similarly.
\end{proof}

\end{document}

%% file: command1.tex
\newcommand{\Sig}{\Sigma}
\newcommand{\ra}{\rightarrow}
\newcommand{\xra}{\xrightarrow}
\newcommand{\ov}{\widehat}

\newcommand{\CFKT}{\text{CFK}}

\newcommand{\B}{\mathbb{B}}

\newcommand{\phibar}{{\overline{\phi}}}
\newcommand{\psibar}{{\overline{\psi}}}

\newcommand{\U}{\mathbb{U}}

\newcommand{\M}{\mathbb{M}}
\newcommand{\E}{\mathbb{E}}

\newcommand\Dual{\mathcal D}
\newcommand\Duality\Dual

\newcommand\RelSpinC{\underline{\SpinC}}
\newcommand\relspinc{s}

\newcommand\x{\mathbf x}

\newcommand\y{\mathbf y}

\newcommand\ModSphere{\ModFlow\left({\mathbb S}\longrightarrow
\Sym^{g-1}(\Sigma_{1})\times \Sym^2(\Sigma_{2})\right)}
\newcommand\ModSpheres\ModSphere

\newcommand\UnparModSp{\widehat \ModSp}
\newcommand\UnparModFlow\UnparModSp

\newcommand\ModMaps{\mathcal M}
\newcommand\ModSp\ModMaps

\newcommand\alphas{\mbox{\boldmath$\alpha$}}

\newcommand\betas{\mbox{\boldmath$\beta$}}

\newcommand\spincrel\relspinc

\newtheorem{thm}{Theorem}[section]

\def\endproof{\relax\ifmmode\expandafter\endproofmath\else
  \unskip\nobreak\hfil\penalty50\hskip.75em\hbox{}\nobreak\hfil\bull
  {\parfillskip=0pt \finalhyphendemerits=0 \bigbreak}\fi}
\def\endproofmath$${\eqno\bull$$\bigbreak}
\def\bull{\vbox{\hrule\hbox{\vrule\kern3pt\vbox{\kern6pt}\kern3pt\vrule}\hrule}}

\newcommand{\Z}{\mathbb{Z}}

\newcommand{\ModSWfour}{\mathcal{M}}
\newcommand{\ModFlow}{\ModSWfour}

\newcommand{\SpinC}{{\mathrm{Spin}}^c}

\newcommand\abuts\Rightarrow
\newcommand\Sym{\mathrm{Sym}}

\newcommand{\etabar}{{\overline{\eta}}}

\newcommand{\Got}{{\phibar^1\otimes \psibar^2}}
\newcommand{\Gto}{{\psibar^1\otimes \phibar^2}}
\newcommand{\Gtt}{{\etabar^1\otimes\etabar^2}}

\newcommand{\Hbb}{{\mathbb{H}}}
\newcommand{\Hinfinf}{{\Hbb^1_\infty\otimes \Hbb^2_\infty}}
\newcommand{\Hinfone}{{\Hbb^1_\infty\otimes \Hbb^2_1}}
\newcommand{\Honeinf}{{\Hbb^1_1\otimes \Hbb^2_\infty}}
\newcommand{\Honeone}{{\Hbb^1_1\otimes \Hbb^2_1}}
\newcommand{\Honezero}{{\Hbb^1_1\otimes \Hbb^2_0}}
\newcommand{\Hzeroone}{{\Hbb^1_0\otimes \Hbb^2_1}}
\newcommand{\Hzerozero}{{\Hbb^1_0\otimes \Hbb^2_0}}

%% file: Note-C-surgery.bbl
\begin{thebibliography}{Dillo 83}

\bibitem[Ef1]{Ef-Whitehead} Eftekhary, E., Longitude Floer homology for
knots and the Whitehead double, \emph{Alg. and Geom.
Topology 5 (2005), also available at} math.GT/0407211
%\bibitem[Ef2]{Ef-gluing} Eftekhary, E., Filtration of Heegaard Floer
%homology and gluing formulas, \emph{preprint, available at}
%math.GT/0410356
\bibitem[Ef2] {Ef-surgery} Eftekhary, E., Heegaard Floer homology and knot surgery,
\emph{preprint, available at} math.GT/0603171

\bibitem[Ef3]{Ef-C-surgery} Eftekhary, E, A combinatorial approach to surgery formulas in Heegaard Floer homology,
\emph{preprint, available at } math.GT/0802.3623
\bibitem[Ef4]{Ef-splicing} Eftekhary, E., Floer homology and splicing knot complements,
\emph{preprint, available at}math.GT/0802.2874
\bibitem[Ef5]{Ef-incompressible} Eftekhary, E., Floer homology and existence of incompressible tori in homology spheres,
\emph{preprint}
\bibitem[Hed]{Hed} Hedden, M., On Floer homology and Berge conjecture on knots admitting lens space surgeries,
\emph{preprint, available at} math.GT/0710.0357v2
%\bibitem[FT]{FT} Fuchs, D., Tabachnikov, S.,
%Invariants of Legandrian and transverse knots in the standard
%contact space,\emph{Topology} vol.36, No.5 (1997) 1025-1053

%\bibitem[Hed]{Hed} Hedden, M., Knot Floer homology of Whitehead doubles,
%\emph{preprint, available at} math.GT/0606094
%\bibitem[Lip]{Lip} Lipshitz, R., Heegaard Floer invariants of
%bordered three-manifolds, \emph{preprint}
\bibitem[MOS]{MOS} Manolescu, C., Ozsv\'ath, P., Sarkar, S.,
A combinatorial description of knot Floer homology,
\emph{preprint, available at} math.GT/0607691
\bibitem[MOST]{MOST} Manolescu, C., Ozsv\'ath, P., Szab\'o, Z., Thurston, D.,
On combinatorial link Floer homology, \emph{preprint, available at}
math.GT/0610559
\bibitem[Ni]{Ni} Ni, Y., Link Floer homology detects the Thurston norm,
\emph{preprint, available at} math.GT/0604360
%\bibitem[OS1]{OS-3mfld}  Ozsv\'{a}th, P., Szab\'o, Z.,
%{Holomorphic disks and topological
%invariants for closed three-manifolds}, \emph{Annals
%of Math.} (2) 159 (2004) no.3 \emph{available at} math.SG/0101206
%\bibitem[OS2]{OS-3m2}  Ozsv\'{a}th, P., Szab\'o, Z.,
%{Holomorphic disks and three-manifold invariants: properties and applications},
%emph{  Annals of Math.} (2) 159 (2004) no.3,
%\emph{available at} math.SG/0105202
\bibitem[OS1]{OS-knot}  Ozsv\'{a}th, P., Szab\'o, Z.,
{Holomorphic disks and knot invariants}, \emph{
Advances in Math.} 189 (2004) no.1, \emph{also available at} math.GT/0209056
\bibitem[OS2]{OS-genus} Ozsv\'{a}th, P., Szab\'o, Z.,
{Holomorphic disks and genus bounds},\emph{ Geom. Topol.}
8  (2004), 311-334
\bibitem[OS3]{OS-surgery} Ozsv\'ath, P., Szab\'o, Z., Knot Floer homology
and integer surgeries, \emph{preprint, available at} math.GT/0410300
\bibitem[OS4]{OS-Qsurgery}Ozsv\'ath, P., Szab\'o, Z., Knot Floer homology and rational surgeries,
, \emph{preprint, available at} math.GT/0504404
%\bibitem[Per]{Per} Perelman, G., Ricci floew with surgery on three-manifolds,
%\emph{preprint, available at} math.DG/0307245
%\bibitem[Ras]{Ras} Rasmussen, J., {Floer homology and knot complements},
%\emph{Ph.D thesis, Harvard univ., also available at}
%math.GT/0306378
\bibitem[SW]{SW} Sarkar, S., Wang, J., A combinatorial description of
some Heegaard Floer homologies, \emph{preprint, available at}
math.GT/0607777
%\bibitem[Sar]{Sar} Sarkar, S., Maslov index of holomorphic triangles,
%\emph{preprint, available at} math.GT/0609673
\end{thebibliography}
